\documentclass[a4paper,reqno]{amsart}

\usepackage{a4wide}
\usepackage{latexsym}
\usepackage[UKenglish]{babel}
\usepackage{graphicx,subfigure}
\usepackage{caption}
\usepackage{algorithm2e}
\usepackage{epstopdf}
	\graphicspath{{./}{figures/}}
\usepackage[colorlinks=true,linkcolor=blue,citecolor=blue]{hyperref}
\usepackage{xcolor}
\usepackage{amsfonts,amsmath,amsthm,mathtools,bbm,cool}
\usepackage{soul} 
\usepackage{enumerate}

\makeatletter
\providecommand{\@LN}[2]{}
\makeatother

\usepackage{tikz}
\usepackage{pgfplots}

\usepackage{hhline}
\usepackage{tabularx} 
\usepackage{multirow}

\usepackage{makecell}
\usepackage{colortbl} 

\newtheorem{theorem}{Theorem}[section]
\newtheorem{proposition}[theorem]{Proposition}
\newtheorem{corollary}[theorem]{Corollary}
\newtheorem{lemma}[theorem]{Lemma}
\newtheorem{remark}[theorem]{Remark}
\newtheorem{assumption}[theorem]{Assumption}
\theoremstyle{definition}

\usepackage{todonotes}

\newcommand{\tbb}{\textcolor{blue}}
\newcommand{\tn}{\textnormal}
\newcommand{\tpp}{\textcolor{purple}}

\newcommand{\numberset}{\mathbb}
\newcommand{\EE}{\numberset{E}} 
\newcommand{\PP}{\numberset{P}}

\newcommand{\B}{\mathcal{B}}
\newcommand{\nuy}{\nu^\mathbf{Y}}

\newcommand{\C}{\mathcal{C}}

\newcommand{\pd}{\partial}

\DeclareRobustCommand{\bbone}{\text{\usefont{U}{bbold}{m}{n}1}} 

\newcommand{\R}{\mathbb{R}}
\newcommand{\Rd}{\mathbb{R}^d}
\newcommand{\Rk}{\mathbb{R}^k}

\usepackage{mathrsfs}

\newcommand\restr[2]{{
		\left.\kern-\nulldelimiterspace 
		#1 
		\littletaller 
		\right|_{#2} 
}}
\newcommand{\littletaller}{\mathchoice{\vphantom{\big|}}{}{}{}}

\newcommand{\xs}{x_{\star}}

\newcommand{\by}{\overrightarrow{\mathbf{y}}}

\newcommand{\bY}{\overrightarrow{\mathbf{Y}}}

\newcommand{\bfy}{\mathbf{y}}
\newcommand{\bfY}{\mathbf{Y}}

\allowdisplaybreaks

\begin{document}
\title{Kinetic variable-sample methods for stochastic optimization problems}

\author{Sabrina Bonandin$^{1*}$}
\address{$^1$Institute for Geometry and Practical Mathematics, RWTH Aachen University, Templergraben 55, Aachen 52062 (Germany)}
\curraddr{}
\email{bonandin@eddy.rwth-aachen.de}
\thanks{The work of SB is supported by the Deutsche Forschungsgemeinschaft (DFG, German Research Foundation) – 320021702/GRK2326 – Energy, Entropy, and Dissipative Dynamics (EDDy).}
\thanks{$^*$\emph{Corresponding author:} Sabrina Bonandin}
\author{Michael Herty$^{2,3}$}
\address{$^2$Institute for Geometry and Practical Mathematics, RWTH Aachen University, Templergraben 55, Aachen 52062 (Germany)\\
$^3$Extraordinary Professor, Department of Mathematics and Applied Mathematics, University of Pretoria, Private Bag X20, Hatfield 0028 (South Africa)}
\curraddr{}
\email{herty@igpm.rwth-aachen.de}
\thanks{MH thanks the Deutsche Forschungsgemeinschaft (DFG, German Research Foundation) for the financial support through 442047500/SFB1481 within the projects B04 (Sparsity fördernde Muster in kinetischen Hierarchien), B05 (Sparsifizierung zeitabhängiger Netzwerkflußprobleme mittels diskreter Optimierung) and B06 (Kinetische Theorie trifft algebraische Systemtheorie).}

\subjclass[2020]{82B40, 65K10, 60K35, 90C26}

\keywords{Global optimization, stochastic optimization problems, particle-based methods, consensus-based optimization, Boltzmann equation, kinetic equations}

\date{\today}

\dedicatory{}

\begin{abstract}
We discuss kinetic-based particle optimization methods and varia\-ble-sample strate\-gies for problems where the cost function represents the expected value of a random mapping. 
Kinetic-based optimization methods rely on a consensus mechanism targeting the global minimizer, and they exploit tools of kinetic theory to establish a rigorous framework for proving convergence to that minimizer. 
Variable-sample strategies replace the expected value by an approximation at each iteration of the optimization algorithm. 
We combine these approaches and introduce a novel algorithm based on instantaneous collisions governed by a linear Boltzmann-type equation.
After proving the convergence of the resulting kinetic method under appropriate parameter constraints, we establish a connection to a recently introduced consensus-based method for solving the random problem in a suitable scaling.
Finally, we showcase its enhanced computational efficiency compared to the aforementioned algorithm and validate the consistency of the proposed modeling approaches through several numerical experiments.
\end{abstract}
\maketitle

\section{Introduction}
\label{sec:intro}

Techniques for solving optimization problems that incorporate uncertain information have become crucial tools in fields such as engineering, business, computer science, and statistics \cite{spall2005introduction,aydin2020structure,strini2019stochasticfinance,lan2020first}.
One approach to formulating such problems is to represent uncertain information using random variables of known probability distribution and considering objective functions involving  quantities such as the expected cost, the probability of violation of some constraint, and variance metrics \cite{schneider2006stochastic}. 
The set of problems resulting from this methodology is commonly known as stochastic optimization problems (SOPs), and, if the optimization effort is undertaken prior to the occurrence of the random event, as static SOPs\footnote{The definition of sSOPs is not unanimous. We refer to \cite{bianchi2009survey}.} (sSOPs). 
In the following, we consider settings in which the objective function involves the expected cost of a random vector $\mathbf{Y}$ defined on the probability space $ (\Omega, \mathcal{A}, \PP) $ and taking values in a set $ E \subset \Rk, k \ge 1$:
\begin{equation}
	\label{eqi: min problem main}
	\min_{x \in \Rd} \{ f(x) := \EE_\PP[F(x,\mathbf{Y})]\}.
\end{equation}
Here, $d \ge 1$,  $F:\Rd \times E \to \R$ is some nonlinear, non-differentiable, non-convex objective function, and, denoting by $\B(E)$ the Borel set of $E$ and by $\nu^{\mathbf{Y}}: \B(E) \to [0,1]$ the law of $\bfY$, $\EE_\PP$ indicates the mathematical expectation with respect to $\PP$, that is, for any $x \in \Rd$,
\[
\EE_\PP[F(x,\mathbf{Y})] = \int_{\Omega} F(x,\mathbf{Y}(\omega)) d\PP(\omega) = \int_E F(x,\mathbf{y}) d\nuy(\mathbf{y}).
\]
We require that $f:\Rd \to \R$ admits a global minimizer $x_{\tn{min}} \in \Rd$.
We assume that for any $x \in \Rd$, $\mathbf{y} \mapsto F(x,\mathbf{y})$ is measurable and $\EE_\PP[F(x,\mathbf{Y})]$ is finite; these are standard assumptions in the context of SOPs (we mention \cite{kall1994stochastic,schneider2006stochastic,birge2011introduction} and the more recent \cite{shapiro2021lectures} for an introduction to the subject). 
We will also refer to \eqref{eqi: min problem main} as the true or original problem.

Typically, the mathematical expectation $f$ in problem \eqref{eqi: min problem main} does not have a closed-form solution, which means that approximation methods for it are required. A common approach is to fix a sample $\by := (\mathbf{y}^{(1)}, \ldots, \mathbf{y}^{(M)}) \in E^M$ of $M \in \mathbb{N}$ realizations of the random vector $\mathbf{Y}$ and $x \in \Rd$, and to approximate $f(x)$ with the Monte Carlo type estimator
\begin{equation}
	\label{defi: fhatM}
	\hat{f}_M(x, \by) := \frac{1}{M} \sum_{j=1}^M F(x, \mathbf{y}^{(j)}),
\end{equation}
also known in the literature as sample average approximation (SAA) \cite{shapiro2003monte, shapiro2021lectures}. $M$ is generally referred to as the sample size.
The main reason this strategy is frequently used lies in the availability of proofs of convergence of optimal solutions and optimal values of the approximated problem $\min_{x\in\Rd} \hat{f}_M(x,\by)$ to those of the true problem \eqref{eqi: min problem main} under fairly general assumptions \cite{shapiro2021lectures,homem2000variable} and through standard probabilistic tools such as the laws of large numbers and the central limit theorem \cite{billingsley2017probability}.
In addition, once $\by$ is fixed, $\hat{f}_M(\cdot,\by)$ is a deterministic function, so the convergence of the strategy is guaranteed by established convergence results of the algorithm used to solve the minimization problem.

The basic concept of the SAA method yields several variations. In this manuscript, we consider variable-sample (VS) techniques \cite{homem2000variable,homem2003variable} in which a new sample $\by_h$ is drawn at each iteration $h$ of the iterative algorithm employed to address the optimization problem. This new sample is used to define $\hat{f}_M$ in the iteration. We remark that this is in contrast to the classical SAA, where a sample is fixed at the beginning, and then the resulting deterministic function is optimized. 
As discussed in \cite{homem2000variable,homem2003variable}, the main feature of VS methods of generating independent estimates of the objective function at different iterations prevents obtaining a candidate minimizer that is strongly dependent on the realization, and, hence, constitutes a clear advantage of VS over SAA methods, where the optimal solutions depend on the initially fixed $\by$. On the other hand, proving the convergence of the former is a challenging task due to the presence of the fluctuating sample $\by_h$, and, thus, due to the function being optimized changing at each iteration.

In recent years, traditional optimization methods have been replaced by meta-heuristics for solving sSOPs \cite{bianchi2009survey,aydin2020structure,juan2023review}. Indeed, although the former are able to find optimal solutions, they are generally only suitable for small problems and require significant computational effort; in contrast, the latter can tackle the complexity and challenges associated with optimization problems under uncertainty, and find good and occasionally optimal solutions \cite{blum2003metaheuristics,bianchi2009survey}. Specifically, meta-heuristics have been combined with SAA and VS methods to solve \eqref{eqi: min problem main}, see for instance \cite{homem2003variable, gutjahr2003converging,gutjahr2004s, homem2000variable, norkin1998optimal} for combinations with the notable meta-heuristics pure random search, ant colony optimization, simulated annealing and branch and bound, respectively.

Among the meta-heuristic-based approaches, the class of consensus-based optimization (CBO) algorithms stands out because of its amenability to a rigorous mathematical convergence analysis \cite{pinnau2017consensus,carrillo2018analytical,carrillo2021consensus,fornasier2021consensus,fornasier2022convergence,totzeck2021trends}. CBO methods consider interacting particle systems that explore the search space $\Rd$ with some degree of randomness while exploiting a consensus mechanism aimed at an estimated minimum. They can be applied to non-smooth functions as a consequence of their derivative-free nature.
Most of the literature available on CBO methods describes and analyses them in two regimes: the finite particle regime (also referred to as microscopic level), see e.g. \cite{ha2020convergence, ha2021convergence, ko2022convergence,bellavia2024noisy}, and the mean-field regime, see e.g. \cite{pinnau2017consensus,carrillo2018analytical,carrillo2021consensus,fornasier2021consensus,fornasier2022convergence}. 
In the former, each particle's time-continuous dynamics is described by a stochastic differential equation (SDE) and its time-discrete counterpart can be derived, for instance, with an explicit Euler-Maruyama scheme \cite{higham2001algorithmic}. The latter is obtained when the number of particles approaches infinity, and provides a statistical description of the interacting system by illustrating how each individual interacts with a theoretical average field generated by all the other particles in the system \cite{pareschi2024optimization}.
Recently, an alternative to the SDEs' formulation at the microscopic level has been proposed in the so-called
kinetic theory-based optimization (KBO) methods \cite{benfenati2022binary}. These methods employ instantaneous binary interaction collisions \cite{villani2002review} between the particles and provide an estimate of the minimizer of the problem by a combination of a local interaction (novelty) and a global alignment (standard consensus mechanism of CBO methods) process. 
As a consequence, the corresponding many-particle dynamics is described by a multidimensional Boltzmann equation  \cite{cercignani2013mathematical}, resulting in a fresh perspective on meta-heuristic optimization \cite{pareschi2024optimization}. This novel approach of using principles of kinetic theory in global optimization has been  applied recently to derive a rigorous converge analysis for the notable meta-heuristics of genetic algorithm \cite{borghi2023kinetic,albi2023kinetic,ferrarese2024localized} and simulated annealing \cite{pareschi2024optimization,borghi2024kinetic}.
Although the approach is innovative, it still connects to the SDEs-based CBO methods at the mean-field level \cite{benfenati2022binary,pareschi2024optimization,borghi2024kinetic}, as the mean-field is derived in the quasi-invariant scaling limit \cite{toscani2006kinetic} inspired by the grazing collisions asymptotics of the Boltzmann equation \cite{desvillettes1992asymptotics,pareschi2013interacting}.

CBO methods are applied to a variety of optimization problems, 
including, e.g., constrained optimization \cite{borghi2023constrained,carrillo2023consensus,fornasier2021consensusspehere,carrillo2024interacting}, multi-objective optimization \cite{borghi2022consensus, borghi2023adaptive}, sampling \cite{carrillo2022consensus}, min-max problems \cite{borghi2024particle,huang2024consensus}, bi-level optimization \cite{trillos2024cb,trillos2024defending}, and problems whose objective is a stochastic estimator at a given point \cite{bellavia2024noisy}.
Recently, they have also been combined with the SAA strategy to address sSOP \eqref{eqi: min problem main} \cite{bonandin2024consensus}. In more detail, a CBO-type algorithm resulting from such combination was proposed, and its consistency with the analogue for the true objective $f$ was recovered at the mean-field level and in the limit of a large sample size. This was achieved by leveraging existing proofs of convergence from SAA theory for the optimal solutions of the approximated problem to those of the true problem.
In the following, we will denote the algorithm resulting from the combination of the SAA strategy and a CBO-type algorithm of \cite{bonandin2024consensus} with \textbf{CBO-FFS}, with FS standing for fixed sample scheme and F standing for fixed sample size. 

To our knowledge, a combination of a consensus-based meta-heuristic
and VS strategies is currently unavailable. In this manuscript, we address this gap by developing the \textbf{LKBO-FVSe} algorithm in Section \ref{sec:LKBO-FVSe}. We take inspiration from the novel approach of the KBO methods of \cite{benfenati2022binary} and propose a microscopic algorithm based on instantaneous collisions. 
We replace the binary interaction formalism of KBO methods by considering the dynamics of the particle system altered through interactions with a scatterer. As a consequence, the particle distribution is modeled using a time-continuous linear Boltzmann equation \cite{case1967linear,villani2002review}.
The weak formulation of the aforementioned equation allows for a theoretical analysis of the model by examining the evolution of observable macroscopic quantities described by a system of ordinary differential equations (ODEs). 
We investigate the stability of the solution involving the first two moments of the particle distribution, namely mean position and energy, and prove convergence to the global minimum $x_{\tn{min}}$ of the true problem, thereby tackling the challenging task of proving convergence of VS-based methods. 
Subsequently, in the quasi-invariant opinion limit, we derive a Fokker-Planck- mean-field-type equation and show its relationship to the mean-field equation derived for CBO-FFS in \cite{bonandin2024consensus}. Finally, having established a relationship between LKBO-FVSe and CBO-FFS, we assess,
through numerical experiments and by exploiting the main advantage of VS over SAA methods of generating independent estimates of the objective function at different iterations, that the former is computationally more efficient than the latter, in the sense specified in Subsection \ref{subsec:fp}.

The rest of the paper is organized as follows. 
In Section \ref{sec:LKBO-FVSe}, we introduce LKBO-FVSe, a variable-sample-inspired algorithm based on instantaneous collisions; we derive the associated Boltzmann-type equation, use it to 
prove that the method is able to capture the global minimum of the true problem, and to draw a connection with the recently introduced CBO-FFS in the quasi-invariant regime. 
In Section \ref{sec:numerics}, we validate the outlined algorithm, show its enhanced efficiency with respect to CBO-FFS, and test the consistency of the proposed modeling approaches through several numerical experiments. We summarize our main conclusions in Section \ref{sec:conclusions} and provide an
overview of possible directions for further research.

\section{A variable-sample kinetic-based algorithm: LKBO-FVSe}
\label{sec:LKBO-FVSe}

In this section, we introduce a consensus-based meta-heuristic that is inspired by variable-sample strategies and whose continuous dynamics is described by a linear kinetic equation of Boltzmann type.
We conduct a comprehensive theoretical analysis at the microscopic (see Subsection \ref{subsec:kineticmodel}), macroscopic (see Subsections \ref{subsec:ev macro}, \ref{subsec:stab}, and \ref{subsec:conv min}), and mean-field level (see Subsection \ref{subsec:fp}).

\subsection{From a microscopic interaction process to a linear Boltzmann model}
\label{subsec:kineticmodel}

Linear Boltz\-mann models arise when a particle system interacts with a scatterer, such as a known distribution of target particles or a fixed obstacle \cite{case1967linear,villani2002review}. Notable examples of models leading to the construction of a linear kinetic equation are the Goldstein-Taylor and radiative transfer ones (see \cite{pareschi2013interacting} and the references therein). 
In the first, two groups of particles traveling along a straight line with velocity $+c$ and $-c$ are considered, and they can simultaneously randomly switch to the opposite velocity. In the second, the charged or neutral particles constituting a particle radiation beam spread out their velocities, for instance uniformly, after the interaction with an obstacle $\mathcal{O} \subset \R^3$.
Both frameworks can be modeled through kinetic considerations, by introducing a random process $V(t)$ describing the proportion of particles with given velocity at time $t$ and by demanding that its variation is due to the instantaneous interaction with a random variable $S$ serving as the external scatterer. In the case of the Goldstein-Taylor model, $S$ is a discrete random variable assuming the values $\pm 1$ with some given probability, 
while, in the case of radiative transfer, $S$ is a variable taking values on $\mathcal{O}$, for instance uniformly distributed. We denote the process resulting from the interaction between $V(t)$ and $S$ as $V'(t)$, and, if $v$, sampled from $V(t)$, is the velocity of a particle before the interaction with the scatterer, the post-interaction velocity is given by
\begin{align*}
	v^{\prime} = sv, \quad &\tn{for the Goldstein-Taylor model}  \\
	v^{\prime} = s \bbone_{\{v \in \mathcal{O}\}} + v \bbone_{\{v \notin \mathcal{O}\}}, \quad &\tn{for the radiative transfer model}
\end{align*}
with $s$ sampled from $S$.

Let us now consider a particle with position $x \in \Rd$ and let $X(t)$ be the random process describing the proportion of particles with position $x$ at time $t$. Let $g(t,x)$ denote its associated probability distribution. 
At this stage, we need to establish an instantaneous rule that generates the position $x^\prime$ from $x$ (or, equivalently, the random process $X'(t)$ from $X(t)$), ensuring that as time progresses, the position converges to the desired minimum. As we are in particular interested in consensus-based meta-heuristics, we expect the dynamics of $x^\prime$ to be governed by the motion of an estimated minimum, usually known as consensus point, which depends on the cost function to be minimized \cite{pinnau2017consensus,benfenati2022binary}, in this case $\hat{f}_M$.
In the setting of variable-sample strategies, a new sample of $M$ entries of $\bfY$ is drawn at each iteration of the optimization procedure according to the distribution $\nu^{\bY}$ \footnote{If, for instance, the $M$ entries are additionally independent, then $\nu^{\bY} = (\nuy)^{\otimes M}$. We remark that all the theoretical results of the manuscript hold for a general sampling distribution $\nu^{\bY}$. In Section \ref{sec:numerics}, we will consider the above product form of $\nu^{\bY}$.} so that in turn the function $\hat{f}_M$ being optimized changes at each iteration. 
Taking inspiration from the Goldstein-Taylor and radiative transfer models, we may see the action of drawing a new sample $\by$ as an interaction with a scatterer $S$ (distributed according to $\nu^{\bY}$), as it leads to changing $\hat{f}_M$, and so $x^\prime$.
Finally, we define the post-interaction position $x^\prime$ as
\begin{equation}
	\label{eq:bin_KBOlinear}
	x' = x + \lambda (x^{\alpha}(t,\by) - x) + \sigma D(t, x,\by) \xi,
\end{equation}
where $x^{\alpha}, \alpha >0,$ is given by
\begin{equation}
	\label{def:xalpha(t,y)}
	x^{\alpha}(t,\by) = \frac{\int_{\Rd} x \omega^{\alpha,\hat{f}_M(\by)}(x) g(t,x) dx}{\int_{\Rd} \omega^{\alpha,\hat{f}_M(\by)}(x) g(t,x) dx}, 
\end{equation} 
with 
\begin{equation}
	\label{def:omegahatfM}
	\omega^{\alpha,\hat{f}_M(\by)}(x) := \exp(-\alpha \hat{f}_M(x,\by)).
\end{equation}
The choice of the exponential weight functions in \eqref{def:xalpha(t,y)} comes from the well-known Laplace principle \cite{dembo2009large} (see also \cite{pinnau2017consensus,carrillo2018analytical,carrillo2021consensus,fornasier2021consensus,fornasier2022convergence}), which guarantees that
$x^{\alpha}(t,\by) \approx \tn{argmin}_{x} \hat{f}_M(x, \by)$ as $\alpha \to + \infty$.
Update rule \eqref{eq:bin_KBOlinear} consists of a balance between a drift term, governed by the constant $\lambda >0$, that pulls the particles toward the temporary estimate of the minimum $x^{\alpha}(t,\by)$, and a diffusion term, driven by the constant $\sigma >0$, that encourages exploration of the search space.  
$\xi$ is a random vector drawn from a normal distribution, and 
$\by$ is a vector sampled from $S \sim \nu^{\bY}$ independently of $\xi$.
$D(\cdot,\cdot,\cdot)$ is a $d \times d$ diagonal matrix characterizing the exploration around the consensus point, that can be of isotropic (all dimensions $l=1,\ldots,d$ are equally explored) or anisotropic type:
\begin{equation*}
	D_{\tn{iso}}(t,x,\by) = |x^{\alpha}(t,\by) -x|I_d, \quad D_{\tn{aniso}}(t,x,\by) = \tn{diag}(x^{\alpha}(t,\by) -x),
\end{equation*}
with $I_d$ denoting the $d$-dimensional identity matrix and $\tn{diag}: \Rd \to \R^{d\times d}$ the operator mapping a vector onto a diagonal matrix with the vector as its diagonal.\\

We hereby present a formal derivation of the Boltzmann-type equation associated to the microscopic collisions \eqref{eq:bin_KBOlinear} and describing the evolution of the particle distribution $g$ over time. We follow the approach taken in \cite{pareschi2013interacting}.

Rule \eqref{eq:bin_KBOlinear} (which we note can be rewritten in terms of the processes $X(t), S$, and $X'(t)$) specifies how to update the position $x$ when an instantaneous interaction with the scatterer $S$ occurs, namely when a sample $\by$ is drawn from $\nu^{\bY}$ at time $t$. 
In order to establish time-continuity, we require that there is some probability of interacting with $S$ which is directly proportional to a short time interval $\Delta t$.
We model this by introducing the random variable $T_{\eta}$, independent of $X(t), S$, and distributed according to a Bernoulli law with parameter $\Delta t/\eta$, for some $\eta > 0$ satisfying the constraint $\Delta t/\eta \le 1$ ($1/\eta$ serves as a measure of the interaction frequency). Then, the time-variation of $X(t)$ is given by 
\begin{equation*}
	X(t+\Delta t) = (1-T_{\eta}) X(t) + T_{\eta} X'(t),
\end{equation*} 
with the law of $X'(t)$ depending in particular on the laws of $X(t), S$, and $\xi$. For any given $\psi \in C^{\infty}(\Rd)$, if we evaluate the mean value of $\psi(X(t+\Delta t))$ and use the independence of $T_{\eta}$ of $X(t), S$, we ultimately recover that the particle distribution $g$ satisfies the weak formulation 
\begin{equation}
	\label{eq:weakBoltz_KBOlinear}
	\frac{d}{dt} \int_{\Rd} \psi(x) g(t,x) dx = \frac{1}{\eta} 
	\int_{\R^{d} \times E^{M}} \left\langle \psi(x') - \psi(x) \right\rangle  g(t,x) \theta_{\bY}(\by) dx d\by
\end{equation}
in the limit $\Delta t \to 0^+$. Here, $\theta_{\bY}$ is the density associated to $\nu^{\bY}$ and $ \langle \cdot \rangle$ denotes the mathematical expectation with respect to $\xi$. 
We point out that, consistently with what done for the Goldstein-Taylor and radiative transfer models, we additionally assumed in the derivation of \eqref{eq:weakBoltz_KBOlinear} that $X(t)$ and $S$ are independent of each other. As commented in \cite{pareschi2013interacting}, this assumption appears quite natural considering the interpretation of the two stochastic processes as a proportion of particles with a given position and scatterer respectively. 
\eqref{eq:weakBoltz_KBOlinear} is complemented with an initial condition $g_0$ satisfying $\int_{\Rd} g_0(x) dx = 1$ and \begin{equation}
	\label{eq:incond_KBOlinear}
	\lim_{t \to 0} \int_{\Rd} \psi(x) g(t,x) dx = \frac{1}{\eta} \int_{\Rd \times E^M} \psi(x) g_0(x) \theta_{\bY}(\by) dxd\by.
\end{equation}

We call the method resulting from the microscopic collision \eqref{eq:bin_KBOlinear} and with Boltzmann equation with weak formulation \eqref{eq:weakBoltz_KBOlinear} LKBO-FVSe. 
The choice of the acronym FVSe is consistent with the nomenclature introduced in \cite{homem2003variable}, with F standing, as in CBO-FFS, for fixed sample size ($M$), VS for variable-sample scheme, and e for equal distribution, as we require the sampling distribution $\nu^{\bY}$ to be time-independent. 
The abbreviation LKBO stands for linear KBO, and emphasizes, on the one hand, the close connection between microscopic collsions \eqref{eq:bin_KBOlinear} and the binary collisions of KBO methods (see also Remark \ref{rem:rel_KBO}). On the other, the fact that the modeling is based on an interaction with a scatterer 
leads to a linear Boltzmann equation.

\begin{remark}
	\label{rem:rel_KBO}
	We note that the definition of $x^{\prime}$ in \eqref{eq:bin_KBOlinear} is a simplified version of the one adopted in the KBO methods of \cite{benfenati2022binary}. In those methods, the use of binary collisions allows for the inclusion of a term that served as a locally weighted best between the two particles involved in the collision, along with the consensus point $x^{\alpha}$. 
\end{remark}

\begin{remark}
	\label{rem:variants_CBO-FVSe}
	The original variable-sample scheme proposed in \cite{homem2000variable,homem2003variable} involves the usage of different sample sizes and sampling distributions along the algorithm. 
	As observed in the paper mentioned, considering a so-called ``schedule of sample sizes" $\{M_h\}_h$ enables a reduction in computational complexity. 
	For example, it allows the user to select a small sample at the initial iterations of the algorithm or to let the algorithm automatically decide what a ``good" sample size is based on statistical tests. 
	The algorithm resulting from this extension is indicated by the acronym VVS, with VS standing for variable-sample scheme and V for variable sample size.
	Subsequently, changing the sampling distribution (in our notation, the law of $S$) as the algorithm progresses through the computations permits, for instance, to use sampling methods that reduce the variance of the resulting estimators.\\
	Although a definition of LKBO-VVS and LKBO-FVS, namely of consensus-based meta-heuristics with the above modifications, is straightforward, we leave the exploration of these variants to future work. 
\end{remark}

\subsection{Evolution of the mean position $m$ and variance $V$}
\label{subsec:ev macro}

Adhering to the strategy presented in \cite{pinnau2017consensus,carrillo2018analytical,benfenati2022binary,albi2023kinetic}, convergence to the global minimum is now proven in three steps: firstly, the weak formulation of the Boltzmann equation \eqref{eq:weakBoltz_KBOlinear} is used to derive ODEs describing the evolution of its first two moments. Then, the existence of a global consensus $\tilde{x}$ and the concentration are proven under minimal conditions on the objective function $F$. Finally, it is shown that  $\tilde{x}$ is a good approximation of $x_{\tn{min}}$. We assume a sufficiently regular and integrable solution $g$ to \eqref{eq:weakBoltz_KBOlinear} exists. For the sake of notational simplicity, we set $\eta = 1$. 

We denote the mean position at time $t \ge 0$ by
\begin{equation}
	\label{def:m}
	m(t) := \int_{\Rd} x g(t,x) dx.
\end{equation}
Using $\psi(x) = x$ in the weak formulation \eqref{eq:weakBoltz_KBOlinear}, we obtain
\begin{equation}
	\label{eq:dt(m)}
	\frac{d}{dt} m(t) = \lambda (x^{\alpha}(t) - m(t)),
\end{equation}
with 
\begin{equation}
	\label{def:xalpha}
	x^{\alpha}(t) = \int_{E^M} x^{\alpha}(t,\by) \theta_{\bY}(\by) d\by.
\end{equation}

Then, we introduce the mean energy at time $t \ge 0$ as
\begin{equation}
	E(t) := \int_{\Rd} |x|^2 g(t,x) dx,
\end{equation}
and the mean variance at time $t \ge 0$ as
\begin{equation}
	\label{def:V}
	V(t) := \frac{1}{2} \int_{\Rd} |x-m(t)|^2 g(t,x) dx = \frac{1}{2} \left( E(t) - |m(t)|^2 \right).
\end{equation}
Using $\psi(x) = |x|^2$ in the weak formulation \eqref{eq:weakBoltz_KBOlinear}, we obtain
\begin{align*}
	\left\langle \varphi(x') - \varphi(x) \right\rangle &= \left\langle |x'|^2 - |x|^2 \right\rangle = \left\langle |x'|^2 \right\rangle - |x|^2 \\
	&= \lambda^2 | x^{\alpha}(t,\by) -x|^2 + \sigma^2 \sum_{r=1}^d D^2_{rr}(t,x,\by) + 2 \lambda x \cdot (x^{\alpha}(t,\by) -x), 
\end{align*}
with $D_{rr}$ the diagonal entry of the matrix $D$ and
\begin{equation*}
	\sum_{r=1}^d D^2_{rr}(t,x,\by) = \kappa |x^{\alpha}(t,\by) -x|^2,
\end{equation*}
where
\begin{equation}
	\label{def:kappa}
	\kappa := d \; \tn{for isotropic} \quad \tn{and} \quad \kappa := 1 \; \tn{for anisotropic exploration.}
\end{equation}
Finally, we deduce 
\begin{equation}
	\label{eq:dt(V)}
	\frac{d}{dt} V(t) =  \frac{1}{2} (\lambda^2 + \kappa \sigma^2) \int_{\Rd \times E^M} |x^{\alpha}(t,\by)-x|^2 g(t,x) \theta_{\bY}(\by) dx d\by - 2 \lambda V(t),
\end{equation}
by additionally using that 
\begin{align*}
	&\lambda \int_{\Rd \times E^M} x \cdot (x^{\alpha}(t,\by) -x) g(t,x) \theta_{\bY} (\by) dx d\by - m(t) \cdot \lambda (x^{\alpha}(t)-m(t))\\ &\quad= -\lambda (E(t) - |m(t)|^2) = - 2 \lambda V(t).
\end{align*}

The evolution of the mean position $m$ and variance $V$ for LKBO-FVSe is described by equations \eqref{eq:dt(m)} and \eqref{eq:dt(V)}, respectively.
As in \cite{pinnau2017consensus,carrillo2018analytical,benfenati2022binary,albi2023kinetic}, we introduce a boundedness assumption on the objective $F$.

\begin{assumption}
	\label{ass:boundedness obj}
	For any $\bfy \in E$, there exist constants $\underline{F}(\bfy), \overline{F}(\bfy) \in \R$ such that
	\begin{equation*}
		\underline{F}(\bfy) \le F(x,\bfy) \le \overline{F}(\bfy) \quad \tn{for any $x \in \Rd$}.
	\end{equation*}
\end{assumption}
\noindent It is straightforward to verify that the hypothesis implies that, for any $x \in \Rd$,
\begin{equation}
	\label{eq:bound_hatfM}
	\underline{f}_M(\by) := \frac{1}{M} \sum_{j=1}^M \underline{F}(\bfy^{(j)}) \le \hat{f}_M(x,\by) \le  \frac{1}{M} \sum_{j=1}^M \overline{F}(\bfy^{(j)}) =: \overline{f}_M(\by).
\end{equation}
\noindent Then, we may give an upper bound for the right-hand side of \eqref{eq:dt(V)} depending exclusively on $V(t)$, the constants $\lambda, \sigma, \alpha$, and $\kappa$ (defined in \eqref{def:kappa}), and the bounds on the cost function $\hat{f}_M$.
\begin{proposition}
	\label{prop:upper dt(V)}
	Let the particle distribution $g(t,x)$ be a weak solution to \eqref{eq:weakBoltz_KBOlinear} with microscopic interactions given by \eqref{eq:bin_KBOlinear}.
	Let $V(t)$ defined in \eqref{def:V} be its variance.
	If the objective function $F$ fulfills Assumption \ref{ass:boundedness obj}, then
	\begin{equation}
		\label{eq:upper dt(V)}
		\frac{d}{dt} V(t) \le - \left(2\lambda - 2(\lambda^2 + \kappa \sigma^2) C_{\alpha} \right) V(t),
	\end{equation}
	for all $t > 0$ and for
	\begin{equation}
		\label{def:Calpha}
		C_{\alpha} = \int_{E^M} e^{\alpha(\overline{f}_M(\by) -\underline{f}_M(\by))} \theta_{\bY}(\by) d\by >0,
	\end{equation}
	with $\underline{f}_M$ and $\overline{f}_M$ defined in \eqref{eq:bound_hatfM}.
\end{proposition}
\noindent The proof of the proposition is based on the following lemma.
\begin{lemma}
	\label{lemma:quadratic upper V}
	If $F$ fulfills Assumption \ref{ass:boundedness obj}, then
	\begin{equation}
		\label{eq:quadratic upper V}
		\int_{\Rd \times E^M} |x^{\alpha}(t,\by)-x|^2 g(t,x) \theta_{\bY}(\by) dx d\by \le 4 C_{\alpha} V(t),
	\end{equation}
	for all $t > 0$ and $C_{\alpha}$ defined in \eqref{def:Calpha}.
\end{lemma}
\begin{proof}[Proof of Lemma \ref{lemma:quadratic upper V}]
	Denote the quadratic term on the left-hand side of \eqref{eq:quadratic upper V} by $Q^{\alpha}(t)$. 
	Using the definition of $x^{\alpha}(t,\by)$ \eqref{def:xalpha(t,y)} and Jensen's inequality, we get
	\[
	Q^{\alpha}(t) \le \int_{\R^{2d} \times E^M} |\xs-x|^2 \frac{\omega^{\alpha,\hat{f}_M(\by)}(\xs)}{|| \:\omega^{\alpha,\hat{f}_M(\by)}(\cdot) \:||_{L^1(g(t,\cdot))}} g(t,\xs) g(t,x) \theta_{\bY}(\by) d\xs dx d\by.
	\]
	It is easy to see that Assumption \ref{ass:boundedness obj}, and in particular condition \eqref{eq:bound_hatfM}, implies that 
	\begin{equation}
		\label{eq:estimate inner xalphaty}
		\frac{\omega^{\alpha,\hat{f}_M(\by)}(\xs)}{|| \:\omega^{\alpha,\hat{f}_M(\by)}(\cdot) \:||_{L^1(g(t,\cdot))}} \le e^{\alpha(\overline{f}_M(\by) -\underline{f}_M(\by))}.
	\end{equation}
	Plugging \eqref{eq:estimate inner xalphaty} in the previous above inequality, we get 
	\begin{align*}
		Q^{\alpha}(t)
		&\le \int_{\R^{2d} \times E^M} |\xs-x|^2 e^{\alpha(\overline{f}_M(\by) -\underline{f}_M(\by))} g(t,\xs) g(t,x) \theta_{\bY}(\by) d\xs dx d\by\\
		&= \left( \int_{\R^{2d}}\!\! |\xs\!-\!x|^2  g(t,\xs) g(t,x) d\xs dx \right)\!\! \underbrace{\left(  \int_{E^M}\!\! e^{\alpha(\overline{f}_M(\by) -\underline{f}_M(\by))} \theta_{\bY}(\by) d\by \right)}_{C_{\alpha} \; \eqref{def:Calpha}}\\
		&= 2C_{\alpha} \int_{\R^{d}} |\xs|^2 g(t,\xs)   d\xs 
		-2C_{\alpha} \int_{\R^{2d}} \xs \cdot x \; g(t,\xs)g(t,x)  d\xs dx\\
		&= 2C_{\alpha} (E(t) -|m(t)|^2) = 4C_{\alpha}V(t).   \qedhere
	\end{align*}
\end{proof}
\begin{proof}[Proof of Proposition \ref{prop:upper dt(V)}]
	The desired estimate on $dV(t)/dt$ is obtained by applying the upper bound from equation \eqref{eq:quadratic upper V} of Lemma \ref{lemma:quadratic upper V} to equation \eqref{eq:dt(V)}.
\end{proof}

We have now all the ingredients to prove concentration and emergence of consensus. The assessment of $\tilde{x}$ being a good approximation of $x_{\tn{min}}$ is postponed to Subsection \ref{subsec:conv min}.
\begin{corollary}
	\label{cor:concentration_consensus}
	Under the assumptions of Proposition \ref{prop:upper dt(V)}, if, additionally,
	\begin{equation*}
		2\lambda - 2(\lambda^2+\kappa \sigma^2) C_{\alpha} > 0,
	\end{equation*}
	then $V(t) \xrightarrow{} 0$ as $t \to +\infty$. In particular, there exists $\tilde{x} \in \Rd$ for which  $m(t) \xrightarrow{} \tilde{x}, x^{\alpha}(t) \xrightarrow{} \tilde{x}$ as $t \to +\infty$.
\end{corollary}
\begin{proof}[Proof of Corollary \ref{cor:concentration_consensus}]
	The result follows from the proof presented in Theorem 4.1 of \cite{carrillo2018analytical} (and its adaptation to the kinetic setting of Corollary 3.1 of \cite{benfenati2022binary}) and Proposition \ref{prop:upper dt(V)}. 
	%
\end{proof}
\begin{remark}
	Note that $C_{\alpha} \xrightarrow{} +\infty$ for $\alpha \to +\infty$, and hence the above condition of Corollary \ref{cor:concentration_consensus} may become unfeasible. We refer to \cite{carrillo2018analytical,benfenati2022binary} for a detailed discussion.
\end{remark}

\subsection{Stability analysis of the equilibrium $(m,V)=(\tilde{x},0)$}
\label{subsec:stab}

We consider system \eqref{eq:dt(m)}-\eqref{eq:dt(V)} and investigate the stability of $(m,V)=(\tilde{x},0)$.

The dependence on $m(t)$ on the right-hand side of $dV(t)/dt$ is seen in the following equivalent system:
\[
\begin{cases}
	\frac{d}{dt} m(t) &= \lambda (-m(t)+ x^{\alpha}(t)),\\
	\frac{d}{dt} V(t) &= ((\lambda^2 + \kappa \sigma^2) -2\lambda)V(t)  + \frac{(\lambda^2 + \kappa \sigma^2)}{2} |m(t)|^2
	- (\lambda^2 + \kappa \sigma^2) x^{\alpha}(t) \cdot m(t)\\
	& \hspace{0.4cm}  + \frac{(\lambda^2 + \kappa \sigma^2)}{2} \int_{E^M} |x^{\alpha}(t,\by)|^2 \theta_{\bY}(\by) d\by.
\end{cases}
\]
The presence of $x^{\alpha}(t)$ in the equation for $dm(t)/dt$ and of the last three terms of the right hand-side of $dV(t)/dt$ makes the ODEs system non-linear and non-autonomous. Therefore, 
we consider the approximated system
\begin{equation}
	\label{eq:system_stability}
	\begin{cases}
		\frac{d}{dt} m(t) &= -\lambda m(t) + \lambda \tilde{x},\\
		\frac{d}{dt} V(t) 
		&= ((\lambda^2 + \kappa \sigma^2) -2\lambda)V(t) + \frac{(\lambda^2 + \kappa \sigma^2)}{2} |m(t)|^2 - \frac{(\lambda^2 + \kappa \sigma^2)}{2} |\tilde{x}|^2.
	\end{cases}
\end{equation}
The above approximation is justified by the following observations.  
For large times, $x^{\alpha}(t)$ converges to $\tilde{x}$ thanks to Corollary \ref{cor:concentration_consensus}. 
In particular, $x^{\alpha}(t) \cdot m(t) \xrightarrow{t \to +\infty} |\tilde{x}|^2$.
\noindent Lastly, 
\begin{align*}
	&\hspace{-0.1cm} \int_{E^M} |x^{\alpha}(t,\by)|^2 \theta_{\bY}(\by) d\by \\ &=\underbrace{\int_{\Rd \times E^M} |x^{\alpha}(t,\by)-x|^2 g(t,x) \theta_{\bY}(\by) dx d\by}_{\tn{$\to 0$ (Lemma \ref{lemma:quadratic upper V})}} + \underbrace{3E(t)}_{\makecell{\tn{\scriptsize{$\to 3|\tilde{x}|^2$ (relation}}\\\tn{\scriptsize{between $E(t)$ and $V(t)$)}}}}\\
	&\quad- \underbrace{2 \; x^{\alpha}(t) \cdot m(t)}_{\tn{$\to 2|\tilde{x}|^2$ (previous step)}}
	\xrightarrow{t \to +\infty} |\tilde{x}|^2
\end{align*}
\noindent We point out that we also numerically verify the consistency of the replacement of the system \eqref{eq:dt(m)}-\eqref{eq:dt(V)} with \eqref{eq:system_stability} in Subsection \ref{subsec:numerics_mom}.

We have as initial conditions
\begin{equation*}
	m(0) = \int_{\Rd} x g_0(x)dx, \quad V(0) = \frac{1}{2} \int_{\Rd} |x-m(0)|^2 g_0(x)dx.
\end{equation*}
Then, if we complement \eqref{eq:system_stability} with the above, the classical theory of existence and uniqueness of solutions to Cauchy problems for ODEs systems (see e.g. \cite{hirsch2013differential,verhulst2012nonlinear}) guarantees that there exists, for all times $t >0$, a unique solution to \eqref{eq:system_stability}: In view of the calculations carried out in Subsection \ref{subsec:ev macro}, the unique solution must be $(m,V)=(\tilde{x},0)$.\\
In order to investigate the stability of such equilibrium, we compute the Jacobian $(d+1) \times (d+1)$ of the right-hand side of \eqref{eq:system_stability} at $(\tilde{x},0)$:
\[
DJ(\tilde{x},0) = 
\begin{pmatrix}
	-\lambda & & 0 & 0\\
	& \ddots & & \vdots \\
	0 & & -\lambda & 0\\
	& (\lambda^2 + \kappa \sigma^2) \tilde{x}^T & & (\lambda^2 + \kappa \sigma^2) - 2\lambda
\end{pmatrix}	
\]
Its eigenvalues are $-\lambda$ with multiplicity $d$ and $(\lambda^2 + \kappa \sigma^2)-2\lambda$ with multiplicity $1$. By construction, $\lambda > 0$. 
Furthermore, $C_{\alpha}$ defined in \eqref{def:Calpha} is by construction greater than $1$, and it holds that $2\lambda - 2(\lambda^2 + \kappa \sigma^2) C_{\alpha} \le 2\lambda - (\lambda^2 + \kappa \sigma^2)$, so the condition of Corollary \ref{cor:concentration_consensus} implies that $(\lambda^2 + \kappa \sigma^2)-2\lambda <0$. Then, all eigenvalues of $DJ(\tilde{x},0)$ are strictly negative.
\begin{corollary}
	\label{cor:stab_result}
	Under the assumptions of Corollary \ref{cor:concentration_consensus}, $(\tilde{x},0)$ is the unique asymptotically stable equilibrium to \eqref{eq:system_stability} with initial conditions $(m(0),V(0))$. In other words, 
	\begin{itemize}
		\item for every neighborhood $\mathcal{O}$ of $(\tilde{x},0)$, there is a neighborhood $\mathcal{O}_1$ of $(\tilde{x},0)$ in $\mathcal{O}$ such that every solution $(m(t),V(t))$ with with initial conditions $(m(0),V(0))$ in $\mathcal{O}_1$ is defined and remains in $\mathcal{O}$ for all $t > 0$, and 
		\item it holds that $\lim_{t \to \infty} (m(t),V(t)) = (\tilde{x},0)$.
	\end{itemize}
\end{corollary}

\subsection{Convergence to the global minimum}
\label{subsec:conv min}

We prove that the global consensus $\tilde{x}$ lies in a neighborhood of the global minimizer $x_{\tn{min}}$ of the true objective $f(\cdot) = \EE_\PP[F(\cdot,\bfY)]$ for appropriately chosen parameters. We again follow the strategy illustrated in \cite{carrillo2018analytical,benfenati2022binary,albi2023kinetic}, and present our modifications in Corollary \ref{cor:convergence}. 

In accordance to the aforementioned papers, we introduce an additional regularity assumption on the objective $F$. We denote by $\C^2(A,B)$ the space of twice continuously differentiable functions from the open set $A \subset \Rd$ to $\R$.
\begin{assumption}
	\label{ass:c2condition}
	For any $\bfy \in E$, $F(\cdot,\bfy) \in \mathcal{C}^2(\Rd,\R)$, and there exist constants $c_1(\bfy),c_2(\bfy)>0$ such that
	\begin{enumerate}
		\item $\sup_{x \in \Rd} |\nabla_x F(x,\bfy)| \le c_1(\bfy)$;
		\item $\sup_{x \in \Rd} |\nabla_x^2 F(x,\bfy)| \le c_2(\bfy)$, where $\nabla_x^2$ denotes the Hessian matrix computed with respect to $x$.
	\end{enumerate} 
\end{assumption}
\noindent As in Subsection \ref{subsec:ev macro}, it is straightforward to verify that Assumption \ref{ass:boundedness obj} implies that, for any $x \in \Rd$,
\begin{equation}
	\label{eq:bound_f}
	\underline{f} := \int_E \underline{F}(\bfy) d\nuy(\bfy) \le f(x) \le  \int_E \overline{F}(\bfy) d\nuy(\bfy) =: \overline{f},
\end{equation}
and that $f$ satisfies Assumption \ref{ass:c2condition} with the constants
\begin{equation}
	\label{def:c1,c2}
	c_1:= \int_E c_1(\bfy) d\nuy(\bfy), \; c_2:= \int_E c_2(\bfy) d\nuy(\bfy).
\end{equation}
\begin{corollary}
	\label{cor:convergence}
	Let the particle distribution $g(t,x)$ be a weak solution to \eqref{eq:weakBoltz_KBOlinear} with microscopic interactions given by \eqref{eq:bin_KBOlinear}.
	Let $F$ fulfill Assumptions \ref{ass:boundedness obj} and \ref{ass:c2condition}, and $\underline{f}$ and $c_1,c_2$ be defined in \eqref{eq:bound_f} and \eqref{def:c1,c2}, respectively.
	If $\lambda,\sigma,\alpha$, $\kappa$, and the initial condition $g_0$ satisfy the inequalities
	\begin{subequations}
		\begin{equation}
			\label{def:mu}
			\mu := 2\lambda - 2(\lambda^2+\kappa \sigma^2) C_{\alpha} > 0,
		\end{equation}
		\begin{equation}
			\label{def:nu}
			\begin{split}
				&\frac{2}{\mu \; ||\omega^{\alpha,f}(\cdot)||_{L^1(g_0)}} \alpha e^{-\alpha \underline{f}} (2\lambda c_1 \sqrt{C_{\alpha}}+(\lambda^2 + \kappa \sigma^2) c_2 C_{\alpha}) \max \{\sqrt{V(0)},V(0)\}\\ &:=\nu < \frac{1}{2},
			\end{split}
		\end{equation}
	\end{subequations}
	with $||\omega^{\alpha,f}(\cdot)||_{L^1(g_0)} := \int_{\Rd} e^{-\alpha f(x)} g_0(x) dx$, 
	then there exists $\tilde{x} \in \Rd$ for which  $m(t) \xrightarrow{} \tilde{x}$ as $t \to +\infty$, and the following estimate on the true objective $f$ holds
	\begin{equation*}
		f(\tilde{x}) \le f(x_{\tn{min}}) + r(\alpha) + \frac{\log 2}{\alpha}
	\end{equation*}
	where $r(\alpha) := -\frac{1}{\alpha} \log ||\omega^{\alpha,f}(\cdot)||_{L^1(g_0)} - f(x_{\tn{min}}) \xrightarrow{\alpha \to \infty} 0$ thanks to the Laplace principle \cite{dembo2009large} if $x_{\tn{min}} \in \tn{supp}(g_0)$. 
\end{corollary}

\begin{proof}
	The first part of the statement is obtained by applying Corollary \ref{cor:concentration_consensus}. 
	For the second part, we follow closely \cite{carrillo2018analytical,albi2023kinetic} and, in particular, Theorem 4.1 of \cite{benfenati2022binary}, 
	with the constants and estimates of Subsection \ref{subsec:ev macro}.
	
	We define $\omega^{\alpha,f}(x) := e^{-\alpha f(x)}$, for some $x \in \Rd$, and 
	\begin{equation}
		\label{def:Malphaf}
		M^{\alpha,f}(t) := \int_{\Rd} \omega^{\alpha,f}(x) g(t,x) dx. 
	\end{equation}
	If we  plug in the choice of test function $\psi(x) = \omega^{\alpha,f}(x)$ in weak formulation \eqref{eq:weakBoltz_KBOlinear}, we get 
	\begin{equation}
		\label{eq:dt(Malpha)}
		\frac{d}{dt} M^{\alpha,f}(t) = \int_{\Rd \times E^M} \left\langle \omega^{\alpha,f}(x') - \omega^{\alpha,f}(x) \right\rangle g(t,x) \theta_{\bY}(\by) dx d\by.
	\end{equation}
	A Taylor expansion of $\omega^{\alpha,f}$ yields
	\begin{align*}
		&\left\langle \omega^{\alpha,f}(x') - \omega^{\alpha,f}(x) \right\rangle\\
		&\quad= 
		\left\langle \nabla_x \omega^{\alpha,f}(x) \cdot (x'-x) + \frac{1}{2} (x'-x) \cdot \nabla_x^2 \omega^{\alpha,f}(\hat{x}) (x'-x) \right\rangle\\
		&\quad \ge -\alpha \lambda e^{-\alpha \underline{f}} c_1 |x^{\alpha}(t,\by)-x| - \frac{\alpha}{2} e^{-\alpha \underline{f}} (\lambda^2 + \kappa \sigma^2) c_2 |x^{\alpha}(t,\by)-x|^2
	\end{align*}
	for $\hat{x} = \gamma x + (1-\gamma)x'$ for some $\gamma \in (0,1)$ and where we have used the bounds \eqref{eq:bound_f} and the boundedness of the second derivatives of $f$ with constants defined in \eqref{def:c1,c2} (we refer to Theorem 4.1 of \cite{benfenati2022binary} for more details on the derivation of the upper estimate).
	Plugging the above bound into \eqref{eq:dt(Malpha)} and using Jensen's inequality, we have
	\begin{align*}
		\frac{d}{dt} M^{\alpha,f}(t) &\ge -\alpha \lambda e^{-\alpha \underline{f}} c_1 \left( \int_{\Rd \times E^M} |x^{\alpha}(t,\by)-x|^2 g(t,x) \theta_{\bY}(\by) dx d\by \right)^{1/2}\\
		&\hspace{0.4cm} -\frac{\alpha}{2} e^{-\alpha \underline{f}} (\lambda^2 + \kappa \sigma^2) c_2 \int_{\Rd \times E^M} |x^{\alpha}(t,\by)-x|^2 g(t,x) \theta_{\bY}(\by) dx d\by.
	\end{align*}
	Thanks to Lemma \ref{lemma:quadratic upper V}, we bound $\int_{\Rd \times E^M} |x^{\alpha}(t,\by)-x|^2 g(t,x) \theta_{\bY}(\by) dx d\by $ by $4 C_{\alpha} V(t)$ and conclude that
	\begin{equation}
		\label{eq:dt(Malpha)up}
		\frac{d}{dt} M^{\alpha,f}(t) \ge -\alpha e^{-\alpha \underline{f}} \left( 2\lambda c_1 \sqrt{C_{\alpha}} + (\lambda^2 + \kappa \sigma^2) c_2 C_{\alpha} \right) \max\{ \sqrt{V(t)}, V(t) \}.
	\end{equation}
	If \eqref{def:mu} holds, we may apply Grönwall's inequality to \eqref{eq:upper dt(V)} and conclude that  $V(t) \le V(0) e^{-\mu t}$ for any $t>0$. 
	Integrating \eqref{eq:dt(Malpha)up} and substituting this conclusion into it, we have
	\begin{align*}
		M^{\alpha,f}(t) &\ge M^{\alpha,f}(0) - \frac{2}{\mu} \alpha e^{-\alpha \underline{f}} \left( 2\lambda c_1 \sqrt{C_{\alpha}} + (\lambda^2 + \kappa \sigma^2) c_2 C_{\alpha} \right) \max\{ \sqrt{V(0)}, V(0) \}\\
		&= M^{\alpha,f}(0) (1-\nu) 
		\underbrace{\ge}_{\eqref{def:nu}} \frac{1}{2} M^{\alpha,f}(0) = \frac{1}{2} ||\omega^{\alpha,f}(\cdot)||_{L^1(g_0)}.
	\end{align*}
	
	Since $m(t) \xrightarrow{} \tilde{x}, V(t) \xrightarrow{} 0$ for large $t$, $M^{\alpha,f}(t) \xrightarrow{t \to +\infty} \omega^{\alpha,f}(\tilde{x})$, and we conclude as in Theorem 4.1 of \cite{benfenati2022binary} by applying Laplace's principle to $f$.
\end{proof}

\subsection{Derivation of a reduced complexity Fokker-Planck model}
\label{subsec:fp}

The weak formulation of the linear Boltzmann equation is employed to calculate the evolution of observable quantities, thereby constructing a bridge between the microscopic interaction \eqref{eq:bin_KBOlinear} and the macroscopic/observable level (see Subsection \ref{subsec:ev macro}). 
Between the two levels lies the so-called mesoscopic/kinetic level, which in our formulation corresponds to the partial differential equation obtained by writing the strong formulation of equation \eqref{eq:weakBoltz_KBOlinear}.
It is well-known that a closed-form analytical derivation of the equilibrium distribution of the kinetic equation is difficult to obtain; for this reason, several asymptotics for it have been proposed to derive reduced complexity models.  
In this subsection, we mention the quasi-invariant opinion limit inspired by the grazing collisions asymptotics of the Boltzmann equation \cite{toscani2006kinetic,pareschi2013interacting}. 
Its key idea is to introduce a scaling parameter that leaves the pre-collisional states unaffected while preserving the model's physical properties. This involves introducing $\varepsilon>0$ and considering the scaling
\begin{equation*}
	t \to \frac{t}{\varepsilon}, \quad  \lambda \to \lambda \varepsilon, \quad \sigma \to \sigma \varepsilon.
\end{equation*} 
Plugging the above scaling into the weak formulation \eqref{eq:weakBoltz_KBOlinear} and letting $\varepsilon \to 0^+$, we get the weak formulation of the Fokker-Planck model
\begin{align*}
	\frac{d}{dt} \int_{\Rd} \psi(x) g(t,x) dx 
	&= \lambda \int_{\R^{d}} \nabla_x \psi(x) \cdot (x^{\alpha}(t)-x)  g(t,x) dx\\
	&+ \frac{\sigma^2}{2} \int_{\Rd} \sum_{r=1}^d \pd^2_{x_r} \psi(x) \left( \int_{E^M} D^2_{rr}(t,x,\by) \theta_{\bY}(\by) d\by \right)  g(t,x) dx,
\end{align*}
with $x^{\alpha}(t)$ given by \eqref{def:xalpha} and $D_{rr}$ denoting the diagonal entry of the matrix $D$.

We have already mentioned in the introduction 
that, over the years, two possible descriptions at the microscopic level have been provided for CBO methods: the SDEs-based one, leading to mean-field-type equations, and the one based on instantaneous collision, leading to kinetic Boltzmann-type equations. We have also mentioned and that the quasi-invariant regime draws a connection between the two. 
In particular, it has already been used in \cite{benfenati2022binary} to assess that the mean-field dynamics of the KBO methods (derived in such limit) is regulated by CBO methods. 
We have also noted that CBO-FFS, an SDEs-based CBO-type method combined with SAA strategies, has been recently introduced to tackle problem \eqref{eqi: min problem main} \cite{bonandin2024consensus}. 
We may now use the reduced complexity model obtained above to create a link between LKBO-FVSe and CBO-FFS. 

For ease of exposition, let us focus on the case of the anisotropic random exploration process (similar considerations can be made for the isotropic case). The strong formulation of the above equation then reads
\begin{equation}
	\label{eq:strongFP_aKBOlinear}
	\begin{split}
		&\pd_t g =  \\
		&\quad \lambda \nabla_x \cdot \left( (x-x^{\alpha}(t)) g \right) + \frac{\sigma^2}{2} \sum_{r=1}^d \pd^2_{x_r} \left( \left( \int_{E^M} (x-x^{\alpha}(t,\by))_r^2 \; \theta_{\bY}(\by) d\by \right) g \right), 
	\end{split}
\end{equation}
with $x^{\alpha}(t,\by)$ given by \eqref{def:xalpha(t,y)}. 
CBO-FFS is obtained by fixing a sample of $M \in \mathbb{N}$ realizations of the random vector $\mathbf{Y}$ $\by = (\mathbf{y}^{(1)}, \ldots, \mathbf{y}^{(M)}) \in E^M$ and by writing a CBO-type algorithm for the deterministic function $\hat{f}_M(\cdot,\by)$ \eqref{defi: fhatM} (we refer to \cite[Section 2]{bonandin2024consensus} for the precise definition of the algorithm). 
In the many-particle limit, it has been shown that the particle distribution 
depends on the fixed sample $\by$  and satisfies the mean-field Fokker-Planck equation
\begin{subequations}
	\label{eq:strongFP_aCBOFFS}
	\begin{equation}
		\label{eq:strongFP_aCBOFFS_eq}
		\pd_t g^{\by} = \lambda \nabla_x \cdot \left( (x-x^{\alpha,g^{\by}}(t,\by)) g^{\by} \right) + \frac{\sigma^2}{2} \sum_{r=1}^d \pd^2_{x_r} \left( (x-x^{\alpha,g^{\by}}(t,\by))_r^2 g^{\by} \right),
	\end{equation}
	\begin{equation}
		\label{def:xalpha(t,y)_gy_CBOFFS}
		x^{\alpha,g^{\by}}(t,\by) = \frac{\int_{\Rd} x \omega^{\alpha,\hat{f}_M(\by)}(x) g^{\by}(t,x) dx}{\int_{\Rd} \omega^{\alpha,\hat{f}_M(\by)}(x) g^{\by}(t,x) dx}
	\end{equation} 
\end{subequations}
(see again \cite[Section 2]{bonandin2024consensus}), where we have made a remark on the dependence of the distribution on the sample fixed by $g^{\by}$. 

If we take a closer look at \eqref{eq:strongFP_aKBOlinear} and \eqref{eq:strongFP_aCBOFFS}, two main differences emerge. First, the particle distribution $g^{\by}$ of CBO-FFS depends on $\by$, while $g$ of LKBO-FVSe is independent of it. Then, both drift and diffusion processes of \eqref{eq:strongFP_aCBOFFS} are driven by a sample-dependent consensus point, while for \eqref{eq:strongFP_aKBOlinear} the processes are driven by $x^{\alpha}(t)$ and a suitable integral of $x^{\alpha}(t,\by)$ in $d\by$, respectively.  
These two differences are expected in light of the main primary distinction between SAA and VS strategies regarding the positioning of the sample drawn. 
Ultimately, the dependence of $g^{\by}$ of CBO-FFS on $\by$ suggests that the output of the method depends strongly on the initially fixed sample $\by$. To address this dependency at the numerical level and thus obtain a candidate minimizer that is independent of $\by$, the authors in \cite{bonandin2024consensus} performed multiple ($n_{\tn{sY}}$) iterations of the algorithm and subsequently averaged the results.  
The two properties on $g$ and the drift and diffusion processes of LKBO-FVSe mentioned above suggest that, conversely, the new algorithm is able to obtain a candidate minimizer that is not dependent on the realization. In particular, we expect that it is not necessary to run the algorithm $n_{\tn{sY}}$ times as required for CBO-FFS: In Subsection \ref{subsec:numerics_ass}, we verify this numerically.

\section{Numerical results}
\label{sec:numerics}

This section is devoted to testing the consistency of the modeling approaches of Section \ref{sec:LKBO-FVSe} by means of several numerical experiments. In more detail, in Subsection \ref{subsec:numerics_val} we show that LKBO-FVSe is able to capture the global minimizer $x_{\tn{min}}$, and, in Subsection \ref{subsec:numerics_ass}, that it is computationally more efficient than CBO-FFS. Then, in Subsection \ref{subsec:numerics_mom}, we show the validity of the approximation done in Subsection \ref{subsec:stab} for several choices of parameters and initial data. 

We first discuss the implementation of LKBO-FVSe based on microscopic collision \eqref{eq:bin_KBOlinear}. 
We simulate it through a direct simulation Monte Carlo (DSMC) method (for a general overview of DMSC methods, we mention e.g. \cite{pareschi2001introduction} and the references therein). We choose the Nanbu scheme for linear Boltzmann equations \cite{nanbu1980direct}.
We fix a time horizon $T>0$ and divide the time interval $[0,T]$ into subintervals of width $\Delta t$ and endpoints $t_h$.
At each time iterate $h$, we consider a collection of $N$ particles $\{x^i_h\}^{i=1,\ldots,N}$ and classify each as interacting (in total $N_c$) or non-interacting. The former updates its position according to collision \eqref{eq:bin_KBOlinear}, while the latter leaves its position unvaried.
%
We summarize the scheme in Algorithm \ref{alg:LKBO-FVSe}.
\RestyleAlgo{ruled}
\begin{algorithm}
	\caption{LKBO-FVSe (microscopic interaction \eqref{eq:bin_KBOlinear})} \label{alg:LKBO-FVSe}
	\textbf{set parameters:} $\lambda, \sigma, \alpha, \Delta t$, $N$, $M$, $\eta$, $\varepsilon$ \;
	\textbf{initialize the positions:} $\{x^i_0\}^{i}$, with $x^i_0 \sim g_0$\;
	$N_c \gets \texttt{Iround}\left( \frac{N \Delta t}{\eta \varepsilon} \right)$\;
	$\tpp{h} \gets 0$; \\
	\While{\tn{stopping criterion on $\tpp{h}$ is not satisfied}}{
		draw one sample: $\by_{\tpp{h}} = (\mathbf{y}^{(1)}_{\tpp{h}},\ldots,\mathbf{y}^{(M)}_{\tpp{h}}) \sim \nu^{\bY}$\;
		compute $\{\hat{f}_{M}(x^i_{\tpp{h}}(\by_{\tn{ag}}),\by_{\tpp{h}}) \}^{i}$ according to \eqref{defi: fhatM}\;
		compute $x^{\alpha}(\tpp{h},\by_{\tpp{h}}) = \sum_{i} x^{\tpp{h}}_h(\by_{\tn{ag}}) \exp(-\alpha \hat{f}_{M}(x^{\tpp{h}}_h(\by_{\tn{ag}}),\by_{\tpp{h}}) )/ \sum_{i} \exp(-\alpha \hat{f}_{M}(x^i_h(\by_{\tn{ag}}),\by_{\tpp{h}}) )$\;
		select $N_c$ colliding particles uniformly among all possible particles\;
		\For{\tpp{$l$} \tn{colliding particle}}{
			sample 	$z^{\tpp{l}}_{\tpp{h}} \sim \mathcal{N}(0,I_d)$ {\footnotesize{(normally distributed with zero mean and identity covariance)}}\;
			update $x^{\tpp{l}}_{\tpp{h}+1}(\by_{\tn{ag}})$: $x^{\tpp{l}}_{\tpp{h}+1}(\by_{\tn{ag}}) = x^{\tpp{l}}_{\tpp{h}}(\by_{\tn{ag}})+\lambda \varepsilon(x^{\alpha}(\tpp{h},\by_{\tpp{h}})- x^{\tpp{l}}_{\tpp{h}}(\by_{\tn{ag}}))  \Delta t + \sigma \sqrt{\varepsilon} D(\tpp{h},x^{\tpp{l}}_{\tpp{h}}(\by_{\tn{ag}}), \by_{\tpp{h}}) \sqrt{\Delta t} z^{\tpp{l}}_{\tpp{h}}$}
		\For{\tpp{$r$} \tn{non-colliding particle}}{update $x^{\tpp{r}}_{\tpp{h}+1}(\by_{\tn{ag}})$: $x^{\tpp{r}}_{\tpp{h}+1}(\by_{\tn{ag}}) = x^{\tpp{r}}_{\tpp{h}}(\by_{\tn{ag}})$}
	}
	\textbf{return:} $ \{x^i_{\tpp{h}}(\by_{\tn{ag}}) \}^{i}_{\tpp{h}}$ 
\end{algorithm}

\noindent 
The algorithm takes as input the scaling factors $\eta > 0$ and $\varepsilon > 0$, introduced in Subsections \ref{subsec:kineticmodel} and \ref{subsec:fp}, respectively. $\eta$ is a measure of the frequency of the sampling action, and, due to its probabilistic interpretation, it must satisfy the constraint $\Delta t/\eta \le 1$; $\varepsilon$ makes the code also suitable for simulating suitable asymptotics of the Bolztmann equation. For the rest of the section, we fix
\begin{equation*} 
	\eta = \Delta t, \quad \varepsilon = 1.
\end{equation*}
The positions are initialized according to $g_{0}$. $\texttt{Iround}(x)$ denotes a suitable integer rounding of a positive real number $x$, and its expression is a consequence of the probabilistic interpretation that underlies the DSMC scheme (we refer to \cite{nanbu1980direct,babovsky1986simulation,pareschi2001introduction} for more details). 
At each iterate $h$, a new sample $\by_h$ is drawn and is used to compute $\hat{f}_M$ and the consensus point. 
We note that the dynamics of the position of particle $i$ depend on all the samples drawn over time (the aggregated sample, denoted by $\by_{\tn{ag}}$ in Algorithm \ref{alg:LKBO-FVSe}), which is why we use the notation $x^i_{h}(\by_{\tn{ag}})$.

As in \cite{bonandin2024consensus}, we choose a cost function $F$ that admits a closed-form expression for the expected value $f$, so that the global minimizer $x_{\tn{min}}$ of $f$ is easily computed: We set
\begin{equation}
\label{def:rastrstoc}
\begin{split}
F(x,\bfY) &= F(x,(Y_1,Y_2)^T)\\ &= \frac{1}{d} \sum_{r=1}^d \left[Y_1(x_r - B)^2 - 10 Y_2\cos (2 \pi (x_r-B)) + 10 \right] + C, 
\end{split}
\end{equation}
and 
\begin{equation}
\label{def:rastr}
\begin{split}
f(x) &= \EE_{\PP} [F(x,\bfY)]\\ &= \frac{1}{d} \sum_{r=1}^d \left[\EE[Y_1](x_r - B)^2 - 10 \EE[Y_2] \cos (2 \pi (x_r-B)) + 10 \right] + C,
\end{split}
\end{equation}
for $x \in \Rd$. 
We observe that if $\EE_{\PP}[Y_1] = \EE_{\PP}[Y_2] = 1$, $f$ coincides with the well-known Rastrigin function with constant shifts $B,C \in \R$ (see e.g. \cite{jamil2013literature}). In the following, we choose $B=C=0$ so that $x_{\tn{min}} = 0 \in \Rd$.

The variable-sample nature of LKBO-FVSe requires to fix two additional parameters compared to the standard CBO algorithms \cite{pinnau2017consensus}: the sample size $M$, and the sampling distribution $\nu^{\bY}$.
We require $\nu^{\bY} = (\nuy)^{\otimes M}$ and assume that $\nuy$ is absolutely continuous with respect to the Lebesgue measure with density $\theta_{\bfY}$. In addition, we impose
\[
\theta_{\bfY}(\bfy) = \hat{\theta}(y_1) \hat{\theta}(y_2), \quad \tn{for any $\bfy = (y_1,y_2)^T \in E \subset \R^2$,}
\]
for some $\hat{\theta}$ probability density function on $\R$. 
Hereafter, we consider
\begin{equation} 
\label{param:CBO-FVSe-tip}
M=50,150,250, \quad \hat{\theta} \sim \mathcal{U}([0.1,1.9]), \mathcal{E}(1), \mathcal{N}(1,1)
\end{equation}
with $\mathcal{U}, \mathcal{E}, and \mathcal{N}$ denoting the uniform, exponential, and normal distributions, respectively. For all three distributions, $\EE_{\PP}[Y_1] = \EE_{\PP}[Y_2] = 1$, and $\hat{\theta}$ is supported on a bounded, semi-infinite, and infinite interval, respectively.   

\subsection{Numerical validation of the kinetic model}
\label{subsec:numerics_val}

We test the performance of LKBO-FVSe on the $20$-dimensional stochastic Rastrigin function \eqref{def:rastrstoc}. We choose a random exploration process of anisotropic type ($D^i_{t, \tn{aniso}}$) as it has been shown to be more competitive than the isotropic one for problems with a high-dimensional search space \cite{carrillo2021consensus,fornasier2022convergence}; we stop the evolution of the algorithm when the final iteration $n_{\tn{it}}$ is reached.

We use two measures for validation, namely the expected success rate, and the error. 
We define a run as successful for LKBO-FVSe if the candidate minimizer 
$x^{\alpha}(n_{\tn{it}},\by_{n_{\tn{it}}})$ is contained in the open $||\cdot||_{\infty}$-ball with radius $\tn{thr}=0.25$
\footnote{The threshold is tuned on the shape of the (stochastic) Rastrigin function in a neighborhood of $x_{\tn{min}}$. See \cite{pinnau2017consensus,carrillo2021consensus} for further details.}
around the true minimizer $x_{\tn{min}}$, and compute the first metric 
by averaging the successful runs over $n_{\tn{CBO}}=100$ realizations of the algorithm. Then, we calculate the expected error 
as the average of $||x^{\alpha}(n_{\tn{it},\by_{n_{\tn{it}}})}-x_{\tn{min}}||_{\infty}$, considering only those runs that have been classified as successful.

We select the parameters for the consensus-based part to be
\begin{equation}
\label{param:test1}
N \!=\! 50, \lambda \!=\!1, \sigma \!=\! 7, \alpha = 30, \; \Delta t = 0.01, n_{\tn{it}} = 10^4, \; d=20, g_0 \sim \mathcal{U}([-3,3]^d).
\end{equation}
We tune them so to achieve a high success rate 
when applying the standard CBO to the Rastrigin function \eqref{def:rastr} \footnote{The definition of a successful run is analogous to that of CBO-FVSe: it suffices to replace the consensus point computed with $\hat{f}_M$ with the one calculated using $f$.}. Specifically, the results we obtain for this case for the two metrics are $98 \%$ and $0.0084$, respectively.

\noindent Then, we choose $M$ and $\hat{\theta}$ according to \eqref{param:CBO-FVSe-tip}, the results for LKBO-FVSe are presented in Table \ref{tab:succ_cbo-FVSe-it10001}. For all choices of $M$ and $\hat{\theta}$, the algorithm yields a high expected success rate and low expected error, with values comparable to the aforementioned those of CBO. We conclude that LKBO-FVSe on $\hat{f}_M$ is able to capture the global minimizer $x_{\tn{min}}$ with a performance similar to that of CBO on $f$.
\begin{table}[h!]
\centering
\begin{tabular}{|c|c|c|c|}\hline
	& $\hat{\theta} \sim \mathcal{U}([0.1,1.9])$ & $\hat{\theta} \sim \mathcal{E}(1)$ & $\hat{\theta} \sim \mathcal{N}(1,1)$ \\ \hline
	$M = 50$ & \cellcolor{yellow!40} $100\%, 0.0081$ $(6145)$  & $98\%, 0.0086$ $(5744)$ & $96\%, 0.0084$ $(6307)$ \\ \hline
	$M = 150$ & $97\%, 0.0081$ $(6159)$ & $99\%, 0.0082$ $(5652)$ & $99\%, 0.0082$ $(6233)$ \\ \hline
	$M = 250$ & $99\%, 0.0082$ $(6050)$ & $98\%, 0.0085$ $(6734)$ & $99\%, 0.0082$ $(6242)$ \\ \hline
\end{tabular}
\caption{Expected success rates and errors at the final iterate $n_{\tn{it}}$ and threshold $0.25$ for LKBO-FVSe with parameters \eqref{param:CBO-FVSe-tip} and \eqref{param:test1} applied to the stochastic Rastrigin function \eqref{def:rastrstoc}. The number in parenthesis represents the iterate at which the success rate of $80 \%$ is reached for the first time (this quantity will be used for Table \ref{tab:succ_cbo-FVSe-it6001}). The purpose of this table is to show the convergence of LKBO-FVSe.}
\label{tab:succ_cbo-FVSe-it10001}
\end{table}

Now we investigate the influence of $M$ and $\hat{\theta}$ on the two metrics. We compute the two  at an iterate for which an intermediate success rate is achieved. As the standard CBO reaches a success rate of $80 \%$ at iterate $5962$, we display the two metrics at the iterate $6000$ in Table \ref{tab:succ_cbo-FVSe-it6001}. We emphasize that Table \ref{tab:succ_cbo-FVSe-it10001} and \ref{tab:succ_cbo-FVSe-it6001} 
share the same initial data $g_0$. As expected, we observe a decrease of the expected success rate and an increase of the expected error as the iteration number decreases (and so by passing from Table \ref{tab:succ_cbo-FVSe-it10001} to \ref{tab:succ_cbo-FVSe-it6001}). 
Subsequently, a close examination of Table \ref{tab:succ_cbo-FVSe-it6001} leads to the following observations. 
If we fix $\hat{\theta}$ and look at the two metrics for varying sample size $M$, we observe that their values remain relatively constant across the rows. This occurs because $M$ enters the algorithm through $\hat{f}_M$, which in turn enters through the consensus point $x^{\alpha}(n_{\tn{it}},\by_{n_{\tn{it}}})$. As a result, variations in the parameter are also affected by the random fluctuations of consensus-type algorithms. Therefore, although $\hat{f}_{250}$ provides a better approximation of $f$ than $\hat{f}_{50}, \hat{f}_{150}$, this is not displayed in the expected success rates and errors.
A similar conclusion can be drawn for fixed $M$ and varying $\hat{\theta}$: The stochastic nature of the algorithm covers the fact that, in the transition from the uniform to the normal distribution, the density support increases, and thus the sample $\by_h$ is drawn from a larger interval.

\begin{table}[h!]
\centering
\begin{tabular}{|c|c|c|c|}\hline
	& $\hat{\theta} \sim \mathcal{U}([0.1,1.9])$ & $\hat{\theta} \sim \mathcal{E}(1)$ & $\hat{\theta} \sim \mathcal{N}(1,1)$ \\ \hline
	$M = 50$ & \cellcolor{yellow!40} $77\%, 0.0084$ & $83\%, 0.0087$ & $76\%, 0.0090$ \\ \hline
	$M = 150$ & $77\%, 0.0086$ & $81\%, 0.0086$ & $78\%, 0.0084$ \\ \hline
	$M = 250$ & $79\%, 0.0091$ & $75\%, 0.0087$ & $74\%, 0.0082$ \\ \hline
\end{tabular}
\caption{Expected success rates and errors at the iterate $6000$ and threshold $0.25$ for LKBO-FVSe with parameters \eqref{param:CBO-FVSe-tip} and \eqref{param:test1} applied to the stochastic Rastrigin function \eqref{def:rastrstoc}. The purpose of this table is to investigate the influence of $M$ and $\hat{\theta}$ on the aforementioned metrics.}
\label{tab:succ_cbo-FVSe-it6001}
\end{table}

\begin{remark}
In this subsection, we have evaluated the expected success rate and error with the threshold equal to $0.25$. We remark that similar values for the metrics may be observed if we consider the lower threshold $0.10$. This observation is further justified by the values attained by the two metrics for the standard CBO with parameters \eqref{param:test1} applied to the Rastrigin function \eqref{def:rastr}: $96\%$ and $0.0079$. 
\end{remark}

\subsection{Increased efficiency of LKBO-FVSe compared to CBO-FFS}
\label{subsec:numerics_ass} 

%
It may be argued that the candidate minimizer 
of LKBO-FVSe and, thus the validation metrics of the previous subsection, are biased by the sample $\by_h$ drawn at each iterate $h$. Indeed, as mentioned in Subsection \ref{subsec:fp}, this was the case for CBO-FFS of \cite{bonandin2024consensus}
and motivated the introduction of a loop on $n_{\tn{sY}}$. In this subsection, we present a variant of LKBO-FVSe (hereafter called LKBO-FVSe-sY), which includes a loop over $n_{\tn{sY}}$. We assess in Table \ref{tab:succ_cbo-FVSe-sY-itvar} that this variant has a similar performance to the one of LKBO-FVSe, hence the algorithm does not require the additional averaging.

We present in Algorithm \ref{alg:LKBO-FVSe-sY} LKBO-FVSe-sY. The idea underlying the new algorithm is simple: We construct the objective $\hat{\hat{f}}_M$ so that it is independent of the sample $\by_h$ and compute the consensus point $x^{\alpha,\hat{\hat{f}}_M}(h)$ with it.
\RestyleAlgo{ruled}
\begin{algorithm}
\caption{LKBO-FVSe-sY} \label{alg:LKBO-FVSe-sY}
\textbf{set parameters:} $\lambda, \sigma, \alpha, \Delta t$, $N$, $M$, $\eta$, $\varepsilon$, $\tbb{n_{\tn{sY}}}$\;
\textbf{initialize the positions:} $\{x^i_0\}^{i}$, with $x^i_0 \sim g_0$\;
$N_c \gets \texttt{Iround}\left( \frac{N \Delta t}{\eta \varepsilon} \right)$\;
$\tpp{h} \gets 0$; \\
\While{\tn{stopping criterion on $\tpp{h}$ is not satisfied}}{
\For{$\tbb{l}_{\tpp{h}}=1,\ldots,\tbb{n_{\tn{sY}}}$}{
	draw one sample: $\by_{\tbb{l}_{\tpp{h}}} = (\mathbf{y}^{(1)}_{\tbb{l}_{\tpp{h}}},\ldots,\mathbf{y}^{(M)}_{\tbb{l}_{\tpp{h}}}) \sim \nu^{\bY}$\;
	compute $\{\hat{f}_{M}(x^i_{\tpp{h}},\by_{\tbb{l}_{\tpp{h}}})\}^i$ according to \eqref{defi: fhatM}\;}
compute $\{\hat{\hat{f}}_M(x^i_{\tpp{h}}) = \frac{1}{\tbb{n_{\tn{sY}}}} \sum_{\tbb{l}_{\tpp{h}}=1}^{\tbb{n_{\tn{sY}}}} \hat{f}_{M}(x^i_{\tpp{h}},\by_{\tbb{l}_{\tpp{h}}})\}^i$\;
compute $x^{\alpha,\hat{\hat{f}}_M}(\tpp{h}) = \sum_{i} x^i_h \exp(-\alpha \hat{\hat{f}}_M(x^i_{\tpp{h}}) )/ \sum_{i} \exp(-\alpha \hat{\hat{f}}_M(x^i_{\tpp{h}}) )$\;
select $N_c$ colliding particles uniformly among all possible particles\;
\For{\tpp{$l$} \tn{colliding particle}}{
	update $x^{\tpp{l}}_{\tpp{h}+1}$ with $x^{\alpha,\hat{\hat{f}}_M}(\tpp{h})$ }
\For{\tpp{$r$} \tn{non-colliding particle}}{update $x^{\tpp{r}}_{\tpp{h}+1}$}
}
\textbf{return:} $ \{x^i_{\tpp{h}} \}^{i}_{\tpp{h}}$ 
\end{algorithm} 

We evaluate the expected success rates and errors for LKBO-FVSe-sY applied to the stochastic Rastrigin function \eqref{def:rastrstoc}, with threshold $0.25$ and parameters \eqref{param:test1} in the second and third columns of Table \ref{tab:succ_cbo-FVSe-sY-itvar}. In view of the conclusions drawn from Table \ref{tab:succ_cbo-FVSe-it6001}, we consider only $M=50$ and $\hat{\theta} \sim \mathcal{U}([0.1,1.9])$, and fix $n_{\tn{sY}} = 50,100$. To facilitate the comparison with the values of the metrics obtained for LKBO-FVSe in Tables \ref{tab:succ_cbo-FVSe-it10001} and \ref{tab:succ_cbo-FVSe-it6001}, we include them in the first column of Table \ref{tab:succ_cbo-FVSe-sY-itvar} and highlight them in yellow. We also consider the iterate $h=8000$ to appreciate how the two metrics evolve during the algorithm's computation.
Ultimately, Table \ref{tab:succ_cbo-FVSe-sY-itvar} was computed alongside the aforementioned tables, and thus shares the same initial data $g_0$.

We observe that the two metrics are fairly stable across the second and third columns, and so, for varying $n_{\tn{sY}}$. This is likely because for a large iterate $h$, the candidate minimizer $x^{\alpha,\hat{\hat{f}}_M}(h)$ is already a good approximation of $x_{\tn{min}}$, which means that $n_{\tn{sY}}$ has minimal influence on its computation. If we compare the results for LKBO-FVSe and LKBO-FVSe-sY, we assess that the metrics perform similarly for all three values of $h$, hence concluding that the loop over $n_{\tn{sY}}$ is not needed for LKBO-FVSe.

\begin{table}[h!]
\centering
\begin{tabular}{|c|c|c|c|}\hline
	& LKBO-FVSe & \makecell{LKBO-FVSe-sY,\\ $n_{\tn{sY}} = 50$} & \makecell{LKBO-FVSe-sY,\\ $n_{\tn{sY}} = 100$}\\ \hline
	$h = n_{\tn{it}}$ & \cellcolor{yellow!40} $100\%, 0.0081$ & $100\%, 0.0084$ & $96\%, 0.0093$ \\ \hline
	$h = 8000$ & $92\%, 0.0086$ & $96\%, 0.0086$ & $94\%, 0.0084$ \\ \hline
	$h = 6000$ & \cellcolor{yellow!40} $77\%, 0.0084$ & $77\%, 0.0084$ & $81\%, 0.0085$  \\ \hline
\end{tabular}
\caption{Expected success rates and errors at three iterates $h$ and threshold $0.25$ for LKBO-FVSe and for LKBO-FVSe-sY with parameters \eqref{param:test1}, $M=50$, and $\hat{\theta} \sim \mathcal{U}([0.1,1.9])$ applied to the stochastic Rastrigin function \eqref{def:rastrstoc}. For LKBO-FVSe-sY, two values of $n_{\tn{sY}}$ were considered. The cells highlighted in yellow are taken from Tables \ref{tab:succ_cbo-FVSe-it10001} and \ref{tab:succ_cbo-FVSe-it6001}.}
\label{tab:succ_cbo-FVSe-sY-itvar}
\end{table}

\subsection{Analysis of consistency of the moments}
\label{subsec:numerics_mom}

We investigate the validity of the replacement of system \eqref{eq:dt(m)}-\eqref{eq:dt(V)} with \eqref{eq:system_stability} 
of Subsection \ref{subsec:stab}. 
We note that both systems are derived for $\eta = 1$. However, if this is not the case, the systems will include an additional scaling factor of $1/\eta$ in front of the right-hand side.
We fix the one-dimensional stochastic Rastrigin function \eqref{def:rastrstoc} ($d=1$).

We address the numerical implementation of the two systems.
We simulate the solution to the true system by replacing the particle distribution $g(t_h,x)$ with the empirical distribution
\begin{equation*}
g^N(t_h,x) = \frac{1}{N} \sum_{i=1}^N \delta\left( x - x^i_h \right) 
\end{equation*}
associated to the collection $\{(x^i_h)\}^{i=1,\ldots,N}$ (computed according to Algorithm \ref{alg:LKBO-FVSe}) in the definition of the mean position $m$ \eqref{def:m} and variance $V$ \eqref{def:V}. This yields
\begin{equation*}
m(t_h) \approx \frac{1}{N} \sum_{i=1}^N x^i_h, \quad V(t_h) \approx \frac{1}{2} \left( \frac{1}{N} \sum_{i=1}^N |x^i_h|^2  -|m(t_h)|^2 \right).
\end{equation*}
We observe that the approximated system depends on the numerical quantity $\tilde{x}$. Thanks to the analysis conducted in Subsection \ref{subsec:conv min}, we know that $\tilde{x} \approx x_{\tn{min}}$, provided that $\alpha$ is sufficiently large, and in the sense of Corollary \ref{cor:convergence}. Then, we may substitute $\tilde{x}$ with $x_{\tn{min}}$ in \eqref{eq:system_stability} and simulate its solution trough the MATLAB ODE solver \texttt{ode45}.

\begin{figure}[h!tbp]
\centering
\subfigure[]{\includegraphics[scale = 0.4]{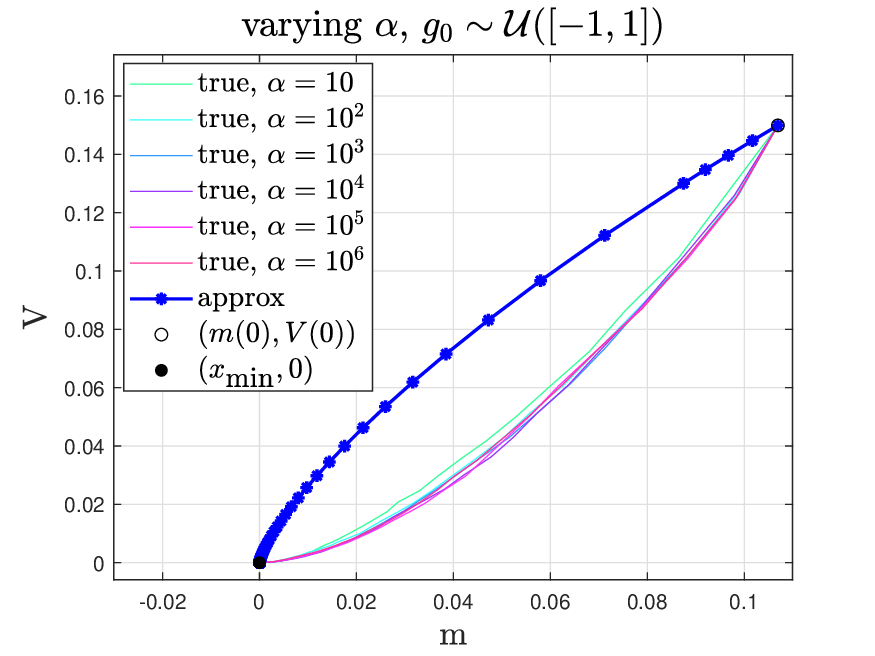}} 
\subfigure[]{\includegraphics[scale = 0.4]{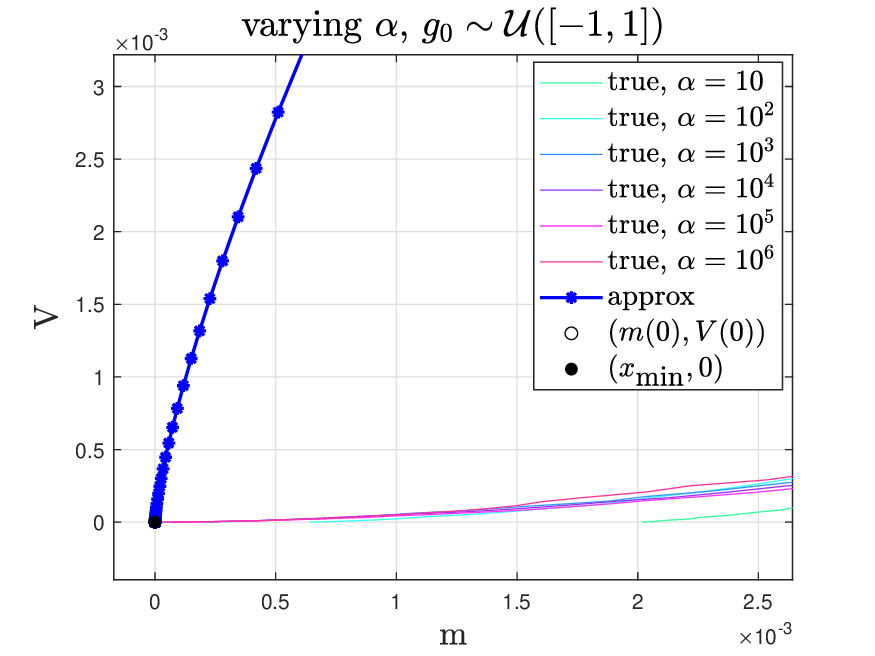}} 
\caption{Plot of the solutions to systems \eqref{eq:dt(m)}-\eqref{eq:dt(V)} and \eqref{eq:system_stability} (with $\tilde{x}$ replaced by $x_{\tn{min}}$) in the phase space $(m,V)$ for several values of $\alpha$. On the right, zoom on a neighborhood of $(x_{\tn{min}},0)$.}
\label{fig:rastrstoc_phfl_varalpha}
\end{figure}

\begin{figure}[h!tbp]
\centering
\includegraphics[scale = 0.45]{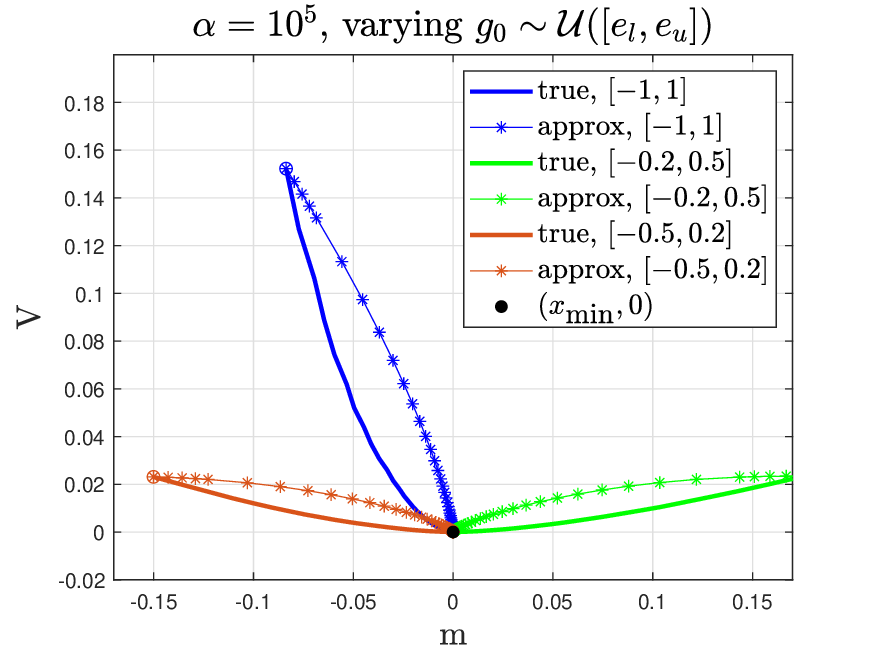}
\caption{Plot of the solutions to systems \eqref{eq:dt(m)}-\eqref{eq:dt(V)} and \eqref{eq:system_stability} (with $\tilde{x}$ replaced by $x_{\tn{min}}$) in the phase space $(m,V)$ for three choices of initial data $g_{0}$.}
\label{fig:rastrstoc_phfl_varin}
\end{figure}

To test the consistency of the approximation of Subsection \ref{subsec:stab}, we plot the solution to the true and approximated systems in the phase space $(m,V)$. We present our results in Figures \ref{fig:rastrstoc_phfl_varalpha} and \ref{fig:rastrstoc_phfl_varin}. The parameters of LKBO-FVSe that are shared between the two figures are
\begin{subequations}
\label{param:test3}
\begin{align}
N = 100; \lambda =&1, \sigma = 0.5, \; \Delta t = 0.1, n_{\tn{it}} = 10^3, \label{param:test3_KBO}\\
&M=50, \hat{\theta} \sim \mathcal{U}([0.1,1.9]).
\end{align}
\end{subequations}

In Figure \ref{fig:rastrstoc_phfl_varalpha}, we fix the initial condition $g_{0} \sim \mathcal{U}([-1,1])$ and consider six values of $\alpha$ (in power of ten) for the true system. We observe that both the trajectories of the true and approximated system decay from the initial datum $(m(0),V(0))$ to the equilibrium $(x_{\tn{min}},0)$ in plot (a). This, in particular, is a graphical representation of the statement of Corollary \ref{cor:stab_result}. In plot (b), we zoom on a neighborhood of $(x_{\tn{min}},0)$ and assess that, the higher the $\alpha$, the closer the end of the trajectory of the true system, namely $(\tilde{x},0)$, to the equilibrium, hence justifying the approximation done in the numerical resolution of system \eqref{eq:system_stability} for $\alpha \gg 1$.

In Figure \ref{fig:rastrstoc_phfl_varin}, we fix $\alpha = 10^5$ 
and choose two other initial conditions in addition to $g_{0} \sim \mathcal{U}([-1,1])$ (identified by the color blue) of the previous figure. More precisely, we fix $g_{0} \sim \mathcal{U}([-0.2,0.5]),$ $\mathcal{U}([-0.5,0.2])$, so we have $(m(0),V(0))$ spanning in both upper quadrants.
Once again, we notice a similar decay in both the true and approximated systems.

\begin{remark}
In Test 3, we couple the use of the stochastic Rastrigin function in one dimension ($d=1$) with parameters \eqref{param:test3_KBO}. A comparison of these parameters with those used for the same objective in dimension $d=20$ from Tests 1 and 2 reveals that the values of $\sigma$ and $n_{\tn{it}}$ are lower. This choice reflects an interplay, already observed in previous works on CBO, between the dimension $d$ of the search space and the parameters $\sigma$ and $n_{\tn{it}}$. 
We refer to \cite{pinnau2017consensus,carrillo2021consensus} for more details concerning the tuning of the parameters.
\end{remark}

\section{Conclusion}
\label{sec:conclusions}

We introduce LKBO-FVSe, a procedure for solving an optimization problem where the cost function is given in the form of an expectation by a combination of variable-sample strategies and kinetic-based meta-heuristics. 
We propose a microscopic description based on instantaneous collisions, with the action of drawing a new sample interpreted as an interaction a scatterer. We model the particle distribution via a time-continuous linear Boltzmann equation, and use the evolution of its moment to prove convergence to the global minimum. We investigate the stability of the solution involving the expected mean and variance, based on an assumption whose consistency we verify numerically. 
We also establish a connection to the recently introduced CBO-FFS for solving the same random problem in the quasi-invariant opinion limit, and show numerically that the feature of LKBO-FVSe of generating independent estimates at each collision leads to a reduction in the computational cost compared to CBO-FFS.

With respect to the standard consensus-based methods, LKBO-FVSe requires fixing two additional parameters, namely the sample size $M$ and the sampling distribution $\nu^{\bY}$. 
In future research, we plan to take up the original idea of variable-sample schemes, and extend the procedure to the usage of a schedule of sample size $\{M_h\}_h$ and varying sampling distribution $\{ \nu^{\bY_h}\}_h$ so as to improve the efficiency of the method proposed.

\section*{Acknowledgments}
SB would like to thank Prof. Claudia Totzeck for her valuable suggestions and insights during the preparation of this manuscript.

\appendix
\bibliographystyle{plain}
\bibliography{references_varkbo}
\end{document}